%% file: ifacconf.tex
\documentclass{ifacconf}

\usepackage{graphicx}      
\usepackage{natbib}        
\usepackage{svg}           
\usepackage{adjustbox}     
\usepackage{amsmath}
\usepackage{amssymb}
\usepackage{enumerate}
\usepackage{bbm}

\begin{document}
\begin{frontmatter}

\title{Hybrid Optimization Methods for Parameter Estimation of Reactive Transport Systems
} 

\thanks[footnoteinfo]{This project has been partly funded by MissionGreenFuels project DynFlex under The Innovation Fund Denmark no. 1150-00001B and the IMI2/EU/EFPIA Joint Undertaking Inno4Vac no. 101007799.}

\author[First]{Marcus Johan Schytt} 
\author[Second]{Halldór Gauti Pétursson} 
\author[First,Third]{John Bagterp Jørgensen}

\address[First]{Department of Applied Mathematics and Computer Science, Technical University of Denmark, Kongens Lyngby, 2800 Denmark (e-mail: \{mschytt, jbjo\}@dtu.dk)}
\address[Second]{MCT Bioseparation ApS, Kongens Lyngby, 2800 Denmark \\(e-mail: hgp@mctbioseparation.com)}
\address[Third]{2-control ApS, Herning, 7400, Denmark}

\begin{abstract}
This paper presents a hybrid optimization methodology for parameter estimation of reactive transport systems. Using reduced-order advection-diffusion-reaction (ADR) models, the computational requirements of global optimization with dynamic PDE constraints are addressed by combining metaheuristics with gradient-based optimizers. A case study in preparative liquid chromatography shows that the method achieves superior computational efficiency compared to traditional multi-start methods, demonstrating the potential of hybrid strategies to advance parameter estimation in large-scale, dynamic chemical engineering applications.
\end{abstract}


\end{frontmatter}

\input{tex/1Introduction}
\input{tex/2Modeling}
\input{tex/3Calibration}
\input{tex/4Chromatography}

\input{tex/5Conclusion}

\begin{ack}
We thank the authors of \cite{chen2024standardized} for providing data and valuable clarifications on the chromatography process. This communication reflects the authors’ views and that neither IMI nor the European Union, EFPIA, or any Associated Partners are responsible for any use that may be made of the information contained therein.
\end{ack}

\bibliography{ifacconf}

\appendix
\input{tex/6Appendix}
\end{document}

%% file: tex/1Introduction.tex
\section{Introduction}
Model-based design strategies have proven highly effective in chemical engineering, enhancing process understanding and enabling rapid optimization~\citep{pistikopoulos2021process}. These strategies rely on mathematical models, often dynamic in nature, to represent the evolving behavior of chemical processes. Chemical processes are generally modeled by differential equations that capture first-principles mass and energy balances, often supplemented with algebraic equations that enforce physical and thermodynamic constraints. To capture a wide variety of systems, models are parametric, with parameters that are often unknown in practice. Unknown parameters must be estimated either from domain knowledge, leveraging insights into parameter ranges and system behavior, or directly from experimental data. When domain knowledge is insufficient or unavailable, parameter estimation can be formulated as an inverse problem. In the case of a dynamic model, this leads to a dynamic optimization problem. Dynamic optimization problems are typically addressed using gradient-based methods for nonlinear programming (NLP). Fortunately, the success of model-based design has been greatly facilitated by the development of efficient and robust NLP solvers~\citep{biegler2018new}.

Accurately calibrating a dynamic model to experimental data is subject to model uncertainties, measurement inaccuracies, parameter correlations, and approximation errors. Together, these factors contribute to a vast and highly irregular optimization landscape with both deep and shallow valleys \citep{arnoud2019benchmarking}. As a result, local (gradient-based) optimization methods, while computationally efficient, tend to capture suboptimal local minima. Their suboptimality often fails to account for the broader model dynamics that can occur across operating conditions in real processes, especially when scaling up from laboratory conditions. Similar challenges arise in process optimization, ultimately limiting its potential for performance improvements. Avoiding suboptimal local minima requires a good initial guess to be fed to the local optimizer, making domain knowledge indispensable. A priori parameter ranges helps narrow the search space, improving the chances of finding a global minimizer. Alternatively, if calibrating a simplified system model is feasible, homotopy methods \citep{schloder1983identification} can be employed to iteratively guide the solution from the simplified system to the original system. In the absence of domain knowledge, a more general approach is necessary.

Global optimization methods address the challenge of local minima by incorporating exploration into the optimization process. However, the key difficulty lies in striking the right balance between exploration and exploitation. Deterministic global optimizers achieve this balance by providing theoretical guarantees for their convergence. Several strategies, such as those based on branch-and-bound or interval arithmetic, have been successfully applied to various dynamic parameter estimation problems in chemical engineering \citep{esposito2000global, singer2006global, lin2006deterministic}. Despite their theoretical guarantees, their main drawback is computational intractability due to the curse of dimensionality. In contrast, stochastic global optimizers, such as population-based metaheuristics, tend to be less affected by this curse. However, they present their own set of challenges. Their performance is highly sensitive to hyperparameters, leading to potentially high computational costs, and tuning these parameters can require considerable effort due to their stochastic nature.

Hybrid optimization methods offer a promising solution to reduce the computational burden of global optimization \citep{arnoud2019benchmarking, villaverde2019benchmarking}. Broadly, a hybrid method may be understood as any method that leverage the strengths of both local and global optimization methods. Global exploration by a global optimizer can identify promising regions of the search space, which are refined through local exploitation by a local optimizer. The most widely applied hybrid methods, multi-start methods, explore the search space through quasi-random sampling (e.g., Sobol sequences, Latin hypercube sampling) to identify candidate solutions for refinement. Despite their popularity, multi-start methods often underperform compared to methods with more strategic search approaches. These employ a global stochastic search phase with strategically timed local refinements to promote faster convergence \citep{egea2007scatter, ugray2007scatter, arnoud2019benchmarking, villaverde2019benchmarking}.

\subsection{Contribution}
In any case, locating global minima of a dynamic optimization problem requires an immense computational effort. This is especially true when processes involve coupled reaction and transport phenomena, where mathematical models expand into systems of partial differential (algebraic) equations (PDEs or PDAEs). High-fidelity models of reactive transport systems account for multiple spatial dimensions, resulting in large-scale discretizations that can limit their practicality in model-based design. To address this, reduced-order models are commonly used to simplify the physics. Among these, one-dimensional advection-diffusion-reaction (ADR) systems offer a powerful framework for various unit operations in chemical engineering, including fixed-bed adsorption columns, catalytic reactors, and membrane separation processes \citep{stegehake2019modeling, kancherla2021modeling}. ADR systems can be efficiently simulated using the method of lines by combining a high-order spatial semi-discretization with a performant numerical integration routine. This approach makes ADR systems well-suited for model-based design of reactive transport processes. However, the use of hybrid optimization strategies for PDE-constrained optimization is rarely addressed in the literature, leaving a gap in its application. This paper contributes by presenting a framework for hybrid optimization that provides a general but simple implementation.

\subsection{Overview}
The rest of the paper is organized as follows: Section~\ref{sec:modeling} introduces the ADR systems and their method of lines discretization. Section \ref{sec:calibration} presents the parameter estimation problem and describes a hybrid optimization methodology. Section \ref{sec:casestudy} presents a case study in chromatography modeling and its parameter estimation by hybrid methods. Finally, Section \ref{sec:conclusion} concludes the paper.

%% file: tex/2Modeling.tex
\section{ADR systems}\label{sec:modeling}
Reduced-order models of reactive transport processes are often simplified to a one-dimensional domain, such as the length of a reactor, pipe, or column. We study the concentrations of $n_{\mathcal{C}}$ reactive components $\mathcal{C}$ over the axial domain $\Omega = [0, L]$ within the time interval $[0, T]$. The mass balances of the component concentrations ${c} = [c_i]_{i \in \mathcal{C}}$ follow a system of semi-linear ADR equations
\begin{align}\label{eq:ADR}
\partial_t {c} = -\partial_z N({c}) + R({c}),
\end{align}
with transportive flux $N$ and reaction term $R = [R_i]_{i \in \mathcal{C}}$. Linear transport is prescribed by splitting the flux into its advective and diffusive parts
\begin{subequations}\label{eq:ADRN}
\begin{align}
    N({c}) &= N_\text{adv}({c}) + N_\text{diff}({c}),\\
    N_\text{adv}({c}) &= {v} \odot {c},\\
    N_\text{diff}({c}) &= -{D} \odot \partial_z {c},
\end{align}
\end{subequations}
where ${v}$ and ${D}$ denote vectors of component velocities and diffusion coefficients, respectively, and $\odot$ denotes the component-wise vector product. The reaction term can be written in stoichiometric form
\begin{align}\label{Rstochiometry}
    R({c}) = \nu^T r({c}),
\end{align}
where $\nu$ is a constant stoichiometric matrix defining reaction kinetics, and $r$ is a vector of kinetic rates. We impose the classical Danckwerts boundary conditions
\begin{subequations}\label{eq:BCs}
\begin{align}
    N({c})(t,0) &= v{c}_\text{in}(t),\\
    N({c})(t,L) &= v{c}(t,L),
\end{align}    
\end{subequations}
where ${c}_\text{in}$ denotes a vector of inlet concentrations. The initial value problem is completed with an initial condition defined by a concentration profile
\begin{align}\label{eq:IC}
    {c}(0,z)={c}_0(z).
\end{align}

\subsection{Discretization}
We discretize the ADR system using the method of lines. We semi-discretize the system and boundary conditions, obtaining a semi-discrete state vector ${x}$ with $n_{x}$ states. The initial concentration profile is similarly discretized and denoted by ${x}_0$. To minimize the computational burden, we use a high-order spatial discretization, specifically the discontinuous Galerkin finite element method (DG-FEM)~\citep{hesthaven2008nodal}. The semi-discrete ADR system can then be formulated as a system of ordinary differential equations (ODEs)
\begin{align}\label{eq:ADRODE}
    \partial_t {x} = f(t, {x}, p), \quad {x}(0)={x}_0,
\end{align}
where $f$ is defined by the discretization and $p$ represent $n_p$ model parameters. The resulting system \eqref{eq:ADRODE} is integrated in time, efficiently, using an explicit singly-diagonally implicit Runge-Kutta (ESDIRK) method \citep{kennedy2016diagonally}.  ESDIRK methods are designed to be stable even in the presence of the highly stiff kinetics common to chemical models.

%% file: tex/3Calibration.tex
\section{Parameter estimation}\label{sec:calibration}
Experimental data is given by a process description and a dataset of measurements. A measurement can generally be expressed as a function of states and model parameters. To ensure sufficient data for statistical inference and to capture the broader model dynamics, it may be necessary to collect measurements from $n_\text{exp}$ different experiments. Each experiment ${\kappa = 1,\dots, n_\text{exp}}$, is defined by an IVP 
\begin{align}
    \partial_t {x}^\kappa = f^\kappa(t,{x}^\kappa, p), \quad {x}^\kappa(0)={x}^\kappa_0,
\end{align}
where the system may depend on experiment-specific boundary conditions, while the model parameters remain constant across all experiments. Measurements of reactive transport systems are typically taken at the domain outlet, with common examples including temperature, pH, and UV absorbance. Each experimental dataset may consist of $n_{\text{disc}}$ different point measurements at $N_\text{disc}^\kappa$ discrete time points, $t_k^\kappa$. These are generally defined as
\begin{align}
    {y}^\kappa_{\text{disc},k}(p) = h_\text{disc}({c}(t_k^\kappa,L), p).
\end{align}
Another common measurement type involves collecting fractions. Effluent fractions are typically measured over fixed time intervals, yielding average or cumulative component concentrations. These $n_\text{cont}$ measurements can be represented by a sequence of $N_\text{cont}^\kappa$ intervals, $[t^\kappa_{\text{start},k},t^\kappa_{\text{end},k}]$, and the corresponding integrals across each interval
\begin{align}
    {y}^\kappa_{\text{cont},k}(p) = \int_{t^\kappa_{\text{start},k}}^{t^\kappa_{\text{end},k}} h_\text{cont}({c}(t,L), p) ,\mathrm{d}t.
\end{align}

In order to estimate the model parameters directly, we minimize a least-squares objective combining the residual sums of squares (RSS) of both the discrete and continuous measurements. The objective function to be minimized with respect to the model parameters is defined as
\begin{equation} \label{eq:least_squares_objective}
    \phi(p) = \sum_{\kappa=1}^{n_{\text{exp}}}\left( \phi^\kappa_\text{disc}(p)+\phi^\kappa_{\text{cont}}(p)\right),  
\end{equation}
where $\phi^\kappa_\text{disc}$ and $\phi^\kappa_\text{cont}$ denote the RSS of both discrete and continuous measurements in experiment $\kappa$, respectively
\begin{subequations}
\begin{align}
    \phi^\kappa_\text{disc}(p) &= \frac{1}{2}\sum_{k=1}^{N_\text{disc}^\kappa} \left\| r_{\text{disc},k}^\kappa(p) \right\|^2,\\
    \phi^\kappa_\text{cont}(p) &= \frac{1}{2}\sum_{k=1}^{N_\text{cont}^\kappa} \left\| r_{\text{cont},k}^\kappa(p) \right\|^2.
\end{align}
\end{subequations}
The residual vectors $r_{\text{disc},k}^\kappa$ and $r_{\text{cont},k}^\kappa$ are defined as
\begin{subequations}
\begin{align}
    r^\kappa_{\text{disc},k}(p) &= w^\kappa_{\text{disc}} \odot \left( {y}^\kappa_{\text{disc},k}(p) - \hat{y}^\kappa_{\text{disc},k} \right),\\
    r^\kappa_{\text{cont},k}(p) &= w^\kappa_{\text{cont}} \odot \left( {y}^\kappa_{\text{cont},k}(p) - \hat{y}^\kappa_{\text{cont},k} \right),
\end{align}
\end{subequations}
where quantities $\hat{y}^\kappa$ denote observed measurements, and $w^\kappa$ denote normalization weights defined component-wise
\begin{subequations}
\begin{align}
    w_{i}^\kappa &= \frac{1}{\max_{k}\{ \hat{y}^\kappa_{i,k}\}}.
\end{align}
\end{subequations}
This weighted least-squares (WLS) formulation normalizes the residuals to ensure consistent magnitudes across measurement types and experiments. Assuming independent, normally distributed errors with zero mean, the WLS solution is equivalent to a maximum likelihood estimate.

To keep the parameters physically meaningful, we specify upper bounds $\overline{p}$, and lower bounds $\underline{p}$. The WLS problem can then be formulated as a dynamic optimization problem
\begin{subequations}
\begin{gather}
    \min_{{x}^\kappa,p}\phi(p) = \sum_{\kappa=1}^{n_{\text{exp}}}\left( \phi^\kappa_\text{disc}(p)+\phi^\kappa_{\text{cont}}(p)\right),  \\
    \partial_t {x}^\kappa = f^\kappa(t,{x}, p), \quad {x}^\kappa(0)={x}^\kappa_0,\\
    \underline{p} \leq p \leq \overline{p}.
\end{gather}
\end{subequations}

\subsection{Metaheuristic optimization}
Population-based metaheuristics offer a black-box approach to optimization. They perform a global heuristic search of the parameter space by maintaining a random population of individuals. The population is evolved over time through operators like selection, mutation, and crossover, with the best individuals carried forward to the next generation. This process continues iteratively until the algorithm meets a stopping criterion, usually a user-defined maximum number of iterations.
Given the ability to evaluate the objective function at relatively low cost through forward simulation, metaheuristics can provide efficient exploration of the parameter space. However, since they do not exploit local information, good solutions may require an excessive number of objective evaluations.

Without a priori knowledge of parameter ranges, the resulting parameter space may be vast. For models based on chemical processes, parameter bounds may be chosen arbitrarily small or large to accommodate unknown kinetics \citep{egea2007scatter}. Since metaheuristics typically initialize by randomly sampling the search space, logarithmic scaling of the parameters can be valuable to select a diverse population over different orders of magnitude. Consider
\begin{align}
    p_i \leftarrow \frac{\log(p_i)-\log(\underline{p}_i)}{\log(\overline{p}_i)-\log(\underline{p}_i)}.
\end{align}

\subsubsection{Multiobjectivization}
Multi-objective metaheuristics extend traditional population-based approaches to handle optimization problems with conflicting objectives. They focus on finding a diverse set of trade-off solutions, using dominance and distance-based ranking to evolve a population toward a well-distributed approximation of the Pareto front. For the problem of parameter estimation, we adopt multi-objective methods through a scheme known as multiobjectivization. Theoretical and practical studies of multiobjectivization have been reviewed in \cite{segura2016using}. Here it is discussed how the methodology can improve the optimization landscape by potentially removing some local optima, and improve the global search by maintaining a diverse population of solutions. By decomposing the original WLS objective into multiple components, we transform the WLS problem into a multi-objective problem. It can be decomposed over $n_{\text{exp}}$ experiments
\begin{align}
    \phi_\text{exp}(p) = \big\{ \phi^\kappa_\text{disc}(p)+\phi^\kappa_{\text{cont}}(p)\big\}_{\kappa},
\end{align}
or over $n_{\text{disc}} + n_{\text{cont}}$ different measurement types
\begin{align}
    \phi_\text{meas}(p) = \bigg\{ \sum_{\kappa=1}^{n_{\text{exp}}} \phi^\kappa_{\text{disc},i}(p), \sum_{\kappa=1}^{n_{\text{exp}}} \phi^\kappa_{\text{cont},j}(p) \bigg\}_{i,j}.
\end{align}
Here we define the RSS of individual measurement types
\begin{align}
    \phi^\kappa_\text{i}(p) &= \frac{1}{2}\sum_{k=1}^{N^\kappa} \left( r_{k,i}^\kappa(p) \right)^2.
\end{align}
In our initial experiments, we found that measurement-based multiobjectivization performs better than the multiobjectivization based on experiment-based decomposition.

\subsection{Local optimization}
The most efficient local optimization methods are gradient-based. They rely on accurate gradient information to be efficient. In dynamic optimization problems, gradients can be computed accurately and efficiently via sensitivity analysis, using either forward or adjoint methods, depending on the size and structure of the problem \cite{petzold2006sensitivity}. Incorporating second-order information, such as the Hessian, can further improve convergence by providing local quadratic convergence. For least-squares problems, the Hessian can be approximated by the Gauss-Newton method. This approximation is computationally straightforward and works well within Levenberg-Marquardt schemes, providing robust performance even for nonlinear problems.

\subsection{Hybrid optimization}
We consider a hybrid optimization method that combines population-based metaheuristics with local refinement. We propose a four-phase framework that iterates for $n_\text{outer}$ outer iterations.

\subsubsection{Candidate search}
An initial population of size $n_\text{pop}$ is randomly generated and iteratively evolved by means of a metaheuristic, either single-objective or multi-objective, for a total of $n_\text{inner}$ inner iterations.

\subsubsection{Candidate selection}
From the current population, the $n_\text{cand}$ candidates with the lowest objective values are selected. For multi-objective methods, candidates are selected based on the lowest sum of objectives to agree with the original single-objective formulation.

To conserve computational resources, selected candidates are filtered using distance criteria as introduced in \cite{ugray2007scatter}. The first, a diversity filter, ensures candidates are sufficiently distinct from previously accepted ones. Given a candidate solution $p$, it is compared against all previously accepted candidates $\tilde{p}$ by verifying
\begin{align}
    \|p-\tilde{p}\| \geq \beta\Delta,    
\end{align}
where $\beta > 0$ is a distance factor, and $\Delta > 0$ represents a minimum critical distance. Candidates satisfying this condition are considered sufficiently diverse and retained for further processing. The next filter, a proximity filter, ensures that selected candidates remain sufficiently distant from refined candidates, which correspond to local minima. For a refined candidate $p^\ast$ and its associated accepted candidate $\tilde{p}$, the radius of its basin of attraction (BOA) is estimated as
\begin{align}\label{eq:boa}
\Delta(p^\ast) = \|\tilde{p} - p^\ast\|.
\end{align}
Candidates are compared against all previously refined candidates by verifying the condition
\begin{align}
    \|p-p^\ast\| \geq \beta\Delta(p^\ast).
\end{align}
If a candidate passes, it is considered to be outside the BOA of previously located minima, and therefore has the potential to locate new minima. If a candidate fails to pass either filter, its potential is disregarded, and the next most promising individual is evaluated instead.

The distance factor is updated between outer iterations, accepting the fact that highly irregular optimization landscapes can have closely spaced local minima all with small BOAs. After all four phases have passed, the distance factor is updated
\begin{align}
    \beta \leftarrow\rho\beta,
\end{align}
for some factor $\rho>0$. In addition, to simplify the selection of an initial critical distance $\Delta$, normalized distances are considered instead
\begin{align}
    \|p\|^2= \frac{1}{n_p}\sum_{i=1}^{n_p}p_i^2.
\end{align}
\subsubsection{Candidate refinement}
Once a candidate is accepted, it undergoes refinement using a local optimization method. To support the filters, its status is updated as $\tilde{p} \leftarrow p$, and its refinement is defined as
\begin{align}
    p^\ast = \text{minimizer}(\tilde{p}).
\end{align}
Lastly, we estimate the radius of its BOA following~\eqref{eq:boa}.

\subsubsection{Candidate replacement}
Once an accepted candidate has been refined, a natural question arises: should we consider the refinements as passive improvements to the metaheuristic's population, or should we actively guide the metaheuristic by reinserting refined candidates into the population? We will refer to the former as passive refinement and the latter as active refinement.

If the passive refinement strategy is chosen, once a refinement is complete, nothing remains but to bookkeep it. On the other hand, if the active refinement strategy is chosen, the refinement seeks to re-enter the metaheuristic's population. This can be done either by directly replacing its associated candidate or by replacing the best individual in the population it dominates. The latter strategy is motivated by its ability to replace previous suboptimal refinements that could mislead the metaheuristic search. Better refinements may represent more accurate solutions, potentially improving exploration in subsequent outer iterations. This strategy outperformed the former in our initial experiments.

Refinements may converge to the same local minima or provide little new information if its neighborhood represents a flat region of the optimization landscape. While filters to address these challenges are discussed in \cite{egea2007scatter}, they are beyond the scope of this paper. Actively manipulating the population of multi-objective metaheuristics is not considered.

%% file: tex/4Chromatography.tex
\section{Chromatography case study}\label{sec:casestudy}
We follow a case study in preparative liquid chromatography presented in \cite{chen2024standardized}. It considers the separation of an industrial antibody feed using ion-exchange chromatography. Chromatography is performed in a fixed bed column in which a mixture of components suspended in a liquid buffer flows through the column. In this case, the feed is a mixture of three species, A1 and A2, which are acidic impurities, and M, the desired product. The components adsorb with different affinities, resulting in different retention rates that provide the means for separation. An eluent component, usually NaCl, is introduced to interfere with the adsorption and provide further control over the elution phase.

\subsection{Chromatography model}
The original authors outline and follow a detailed approach for accurate model development. We implement the column model described in \cite{chen2024standardized}, modeling seven component concentrations, with minor modifications. Instead of accounting for extra-column effects, we time-shift the measurements to adjust for delay volumes in the real system. In addition, the model isotherm is specified by steric mass action and has twelve parameters that must be estimated from data. Appendix~\ref{app:chromamodel} describes how the chromatography model fits within the framework of the ADR systems presented in Section \ref{sec:modeling}.

\subsection{Process description}
In the original case study, bind-and-elute experiments are performed and measurements of the output UV absorbance and effluent fractions are collected for model calibration. The bind-and-elute experiments are run as a three-phase batch process after the column is initially equilibrated with a low-salt buffer.
\begin{enumerate}
    \item Loading phase: The feed mixture is suspended in the low-salt buffer and injected into the column. 
    \item Wash phase: The low-salt buffer is injected to wash out weakly bound impurities.
    \item Elution phase: A high-salt buffer is gradually introduced to elute the feed components.
\end{enumerate}
Five gradient elution experiments are performed in which the high-salt buffer concentration is increased linearly during the elution phase. In addition, one single-step elution experiment is performed in which the high-salt buffer concentration is stepped at the beginning of the elution phase. The same feed mixture is used for all experiments.

\subsection{Parameter estimation}
We aim to calibrate the chromatography model using the first three gradient elution experiments. They are characterized by low loading volumes and provide valid unsaturated UV absorbance measurements. The last single-step elution experiment is used for validation. The UV measurements are modeled using the Beer-Lambert law
\begin{align}
    y^\kappa_{\text{UV},k} = \sum_{i \in \{\text{A1,A2,M}\}}\gamma_i c_i(t_k^\kappa,L),
\end{align}
with the scaling factor $\gamma_i = M_{w,i} \lambda_i \ell$, defined by the molecular weight $M_{w,i}$, the extinction coefficient $\lambda_i$, and the UV path length $\ell$. Each experimental UV dataset contains several thousand measurements, which we trim and subsample to a concise set of $100$ data points for analysis. In contrast, fraction data is typically scarce, with approximately fifteen fractions per peak. Fractions are modeled by averaging the integral effluent concentration
\begin{align}
    y^\kappa_{\text{frac},k,i} = \frac{1}{t^\kappa_{\text{end},k} - t^\kappa_{\text{start},k}} \int_{t^\kappa_{\text{start},k}}^{t^\kappa_{\text{end},k}} c_{i}(t, z)\,\mathrm{d}t.
\end{align}
The model was calibrated in \cite{chen2024standardized} by a sequential parameter-by-parameter (PbP) estimation technique based on protein retention times. While the PbP method itself requires only simple linear regression, its data requirements are more stringent limiting its general use. In our numerical experiments, we evaluate the ability of both metaheuristics and hybrid methods to calibrate the chromatography model. Using (refined) PbP-calibrated parameters as a baseline, we compare the strategies outlined in Section \ref{sec:calibration}.

\subsection{Numerical experiments}

\begin{table}[tb]
\caption{Baseline objective values.}
\centering
\begin{tabular}{ccc}\label{tb:baseline}
Multi-start & Reference & Refined \\ \hline
21.4299     & 1.7705    & 0.8875   
\end{tabular}%
\end{table}

To balance accuracy and efficiency, we discretize the column into five elements with DG-FEM, using third-order polynomials within each element. This results in a semi-discrete state with 140 variables. Additionally, integration tolerances are set to $\varepsilon_\text{rel}=10^{-3}$ and $\varepsilon_\text{abs}=10^{-7}$, adjusted for problem scaling, and handled by a three-stage second-order ESDIRK method TR-BDF2 \citep{hosea1996analysis}. We use staggered-direct forward sensitivity analysis to generate sensitivities \citep{caracotsios1985sensitivity}. The model and its discretization have been implemented in Julia 1.11, and all calculations were performed serially on a an Intel~Xeon~Gold~6342~(2.80~GHz) workstation. With this setup, a nominal simulation of the three experiments took approximately $150$ milliseconds, while the corresponding sensitivity analysis took around $3$ seconds. To ensure a fair comparison between hybrid and metaheuristic methods, we impose a fixed computational budget of two CPU hours. Accordingly, we allow the metaheuristics to perform $48\,000$ function evaluations, assuming no overhead. If we restrict each hybrid method to $20$ refinements, each with a time limit of $60$ seconds, they are allowed a total of $40\,000$ function evaluations for their search. We choose the evolutionary centers algorithm (ECA) \citep{mejia-de-dios2019new} as our single-objective metaheuristic, while we use the non-dominated sorting genetic algorithm (NSGA-II) \citep{deb2002fast} for measurement-based multiobjectivization. We make no effort in their tuning, simply setting their population size to $n_\text{pop}=200$. As our local optimizer, we use the trust-region Newton solver TRON \citep{lin1999newtons}, with a Gauss-Newton Hessian approximation. We consider both passive and active refinement of the metaheuristics. Their filters use a critical distance of $\Delta = 0.1$ and an initial factor $\beta=1$, updated by $\rho = 10^{-11/n_\text{outer}}$, yielding a minimum critical distance of $10^{-12}$. To fit the computational budget, the number of inner and outer iterations depends on the number of candidates. We fix $n_\text{cand}=2$ for both the active and passive methods. As our primary benchmark, we employ a conventional multi-start method utilizing Sobol sequences. This method, like the other hybrid approaches, is permitted a total of $20$ refinements initiated from the best-performing solutions among $40\,000$ candidate evaluations. Its best objective value is reported in Table \ref{tb:baseline} together with the objective values of the PbP-calibrated parameters from \cite{chen2024standardized} and their gradient-based refinement.

\subsubsection{Results}

\begin{table}[tb]
\caption{Best objective values.\\(P) Passive. (A) Active.}
\centering
\begin{tabular}{ccc|cc}
\label{tb:results}
ECA    & ECA(P) & ECA(A) & NSGA-II & NSGA-II(P) \\ \hline
0.8592 & 0.8605 & 0.8656 & 0.9553  & 0.8967    
\end{tabular}%
\end{table}

\begin{figure}[tb]
\begin{center}
\adjustbox{trim=0 0.25cm 0 0.25cm}{\includegraphics[width=8.4cm]{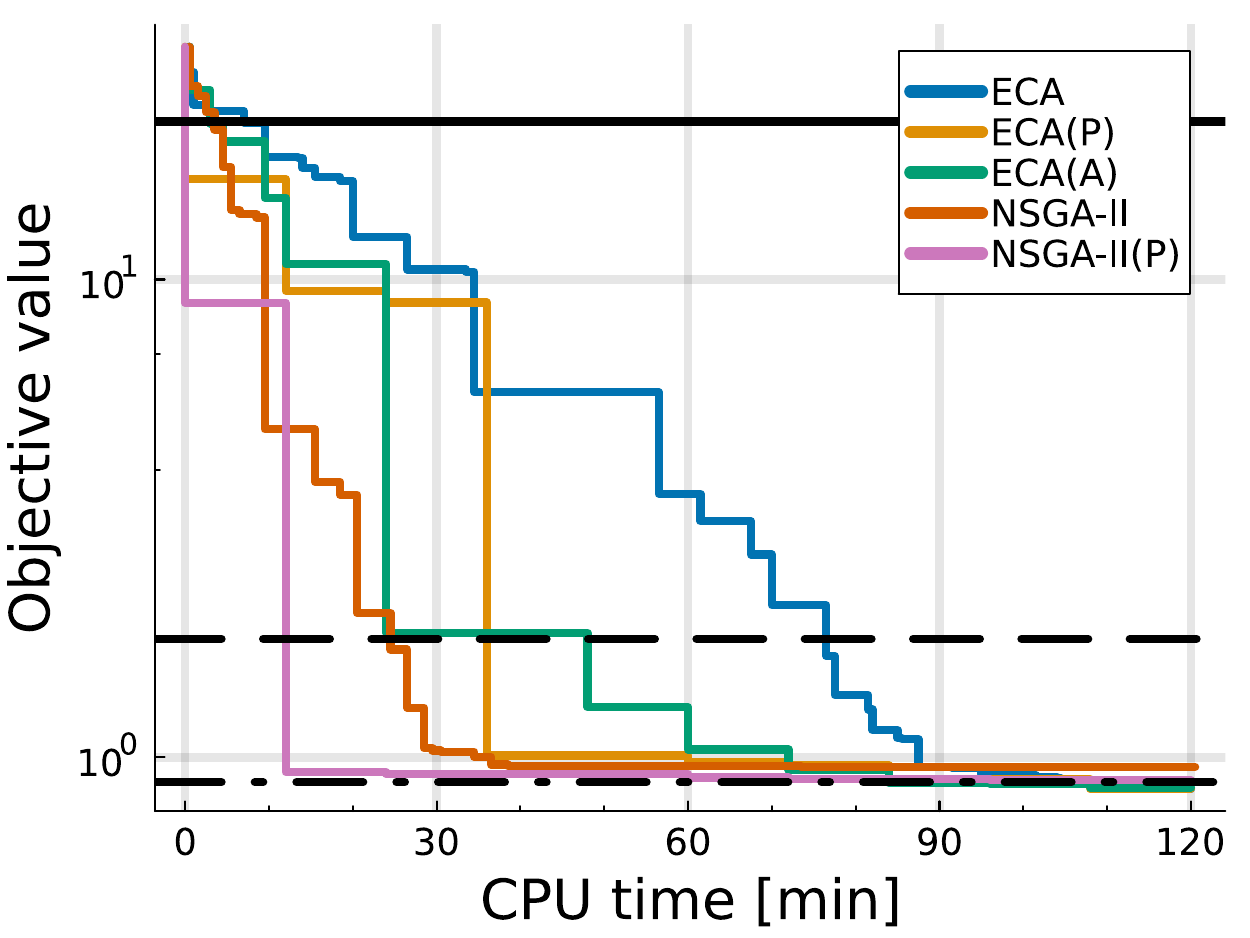}}    
\caption{Optimizer convergence. Multi-start (solid).\\ Reference (dashed). Refined (dot-dashed).} 
\label{fig:results}
\end{center}
\end{figure}

\begin{figure}[tb]
\begin{center}
\adjustbox{trim=0 0.25cm 0 0.25cm}{\includegraphics[width=8.4cm]{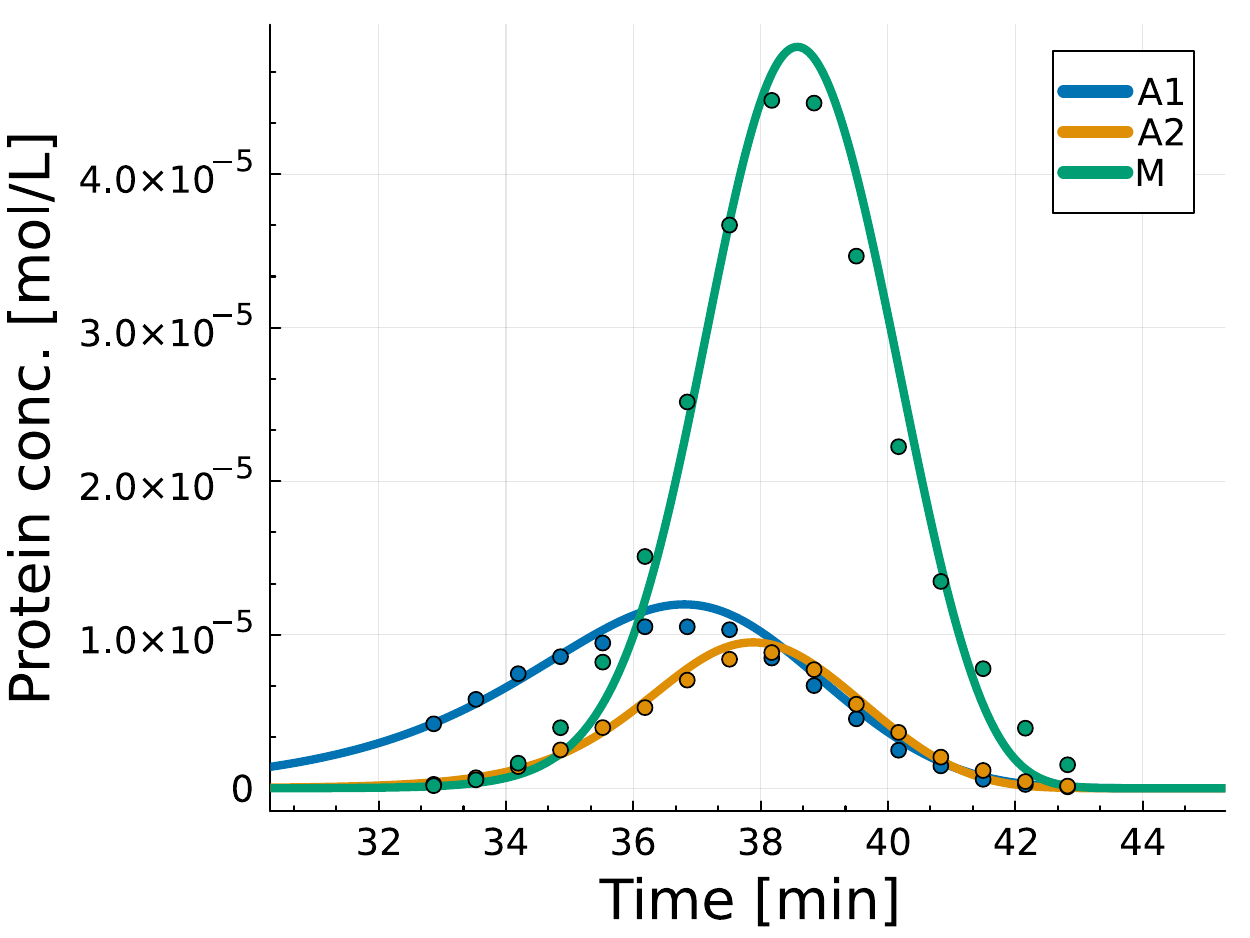}}    
\caption{Simulated outlets against gradient elution experiment fraction data. Simulation (lines). Data (circles).} 
\label{fig:sim}
\end{center}
\end{figure}

\begin{figure}[tb]
\begin{center}
\adjustbox{trim=0 0.25cm 0 0.25cm}{\includegraphics[width=8.4cm]{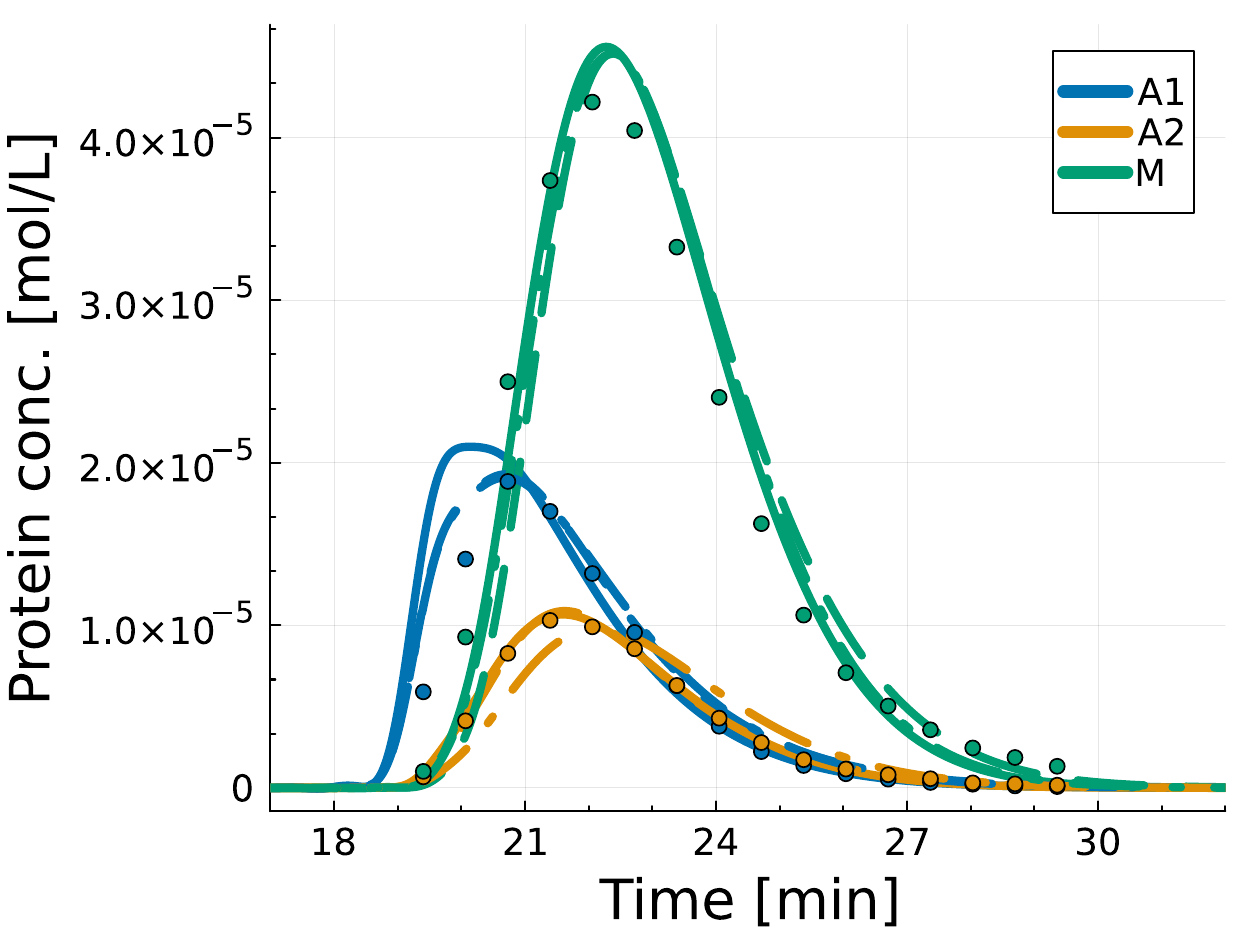}}    
\caption{Simulated outlets against validation data. Refined (solid). ECA(P) (dashed). NSGA-II(P) (dot-dashed).} 
\label{fig:valid}
\end{center}
\end{figure}

The best objective values obtained within the computational budget by each optimizer are shown in Table~\ref{tb:results}. These values illustrate that a simple multi-start method is not sufficient for the problem at hand. In fact, much better performance can be achieved with population-based metaheuristics alone. In terms of objective values, both ECA and NSGA-II manage to beat the reference, but only ECA beats the refined reference. Let us now consider the hybrid passive and active refinement strategies. To get a better sense of the results, we examine their convergence shown in Figure \ref{fig:results}. Including ECA in a hybrid strategy undoubtedly accelerates its initial rate of convergence. However, it does so at the expense of higher objective values. This trade-off may be explained by the reduced exploration performed by hybrid strategies. Indeed, ECA alone only achieves a better objective value after the initial $40\,000$ objective evaluations. It remains inconclusive whether a passive or an active refinement strategy is preferred. Active refinement seems to push the metaheuristic towards faster convergence, while the passive strategy manages to find a better local minimum. On the other hand, when we apply passive refinement to NSGA-II, we see both an acceleration of convergence, and an improvement in objective value. Since multi-objective metaheuristics seek a diverse representation of the Pareto front, we cannot expect NSGA-II alone to sacrifice population diversity for a better objective value. Therefore, its passive refinement consistently finds better objective values. Nevertheless, the refined reference, the metaheuristics, and our proposed hybrid methods all show excellent and identical model fits to the training data. This is shown in Figure \ref{fig:sim}. However, as shown in Figure \ref{fig:valid}, the models differ in their fit to the validation data. Here, the two hybrid methods agree in their fit to the first acidic component A1, while the refined reference and passively refined ECA models agree on the second acidic component A2. All three models predict the desired product M similarly.

%% file: tex/5Conclusion.tex
\section{Conclusion}\label{sec:conclusion}
This paper has introduced a framework for hybrid optimization, and demonstrated its ability to accelerate the parameter estimation of reactive transport systems modeled by ADR systems. A case study in chromatography illustrates that metaheuristics and their acceleration by hybrid optimization methods can be successful in cases where traditional multi-start methods fail. While better results can be achieved through algorithm tuning and code parallelization, considerable potential lies in the use of simple hybrid optimization methods and their integration with multi-objective optimization.

%% file: tex/6Appendix.tex
\section{ADR formulation of chromatography}\label{app:chromamodel}
The equilibrium dispersive model used in \cite{chen2024standardized} fits within the framework of the ADR systems presented in Section \ref{sec:modeling}. The column is
characterized by its length $L$, and is modeled axially $\Omega = [0, L]$. The system includes four components
    $\mathcal{C}_q = \{\text{A1}, \text{A2}, \text{M}\}$,
each with their own mobile and stationary phase. Including a salt, the component mixture $\mathcal{C}$ is described by a vector of mobile and stationary phase concentrations $c_i$ and $q_i$, respectively 
\begin{align}
    c = [c_\text{NaCl}, c_\text{A1}, c_\text{A2}, c_\text{M}, q_\text{NaCl},q_\text{A1}, q_\text{A2}, q_\text{M}].
\end{align}
Their transport is governed by constant advection and diffusion motivated by Darcy's law of fluid flow in porous media and Fick's law of diffusion. This corresponds to \eqref{eq:ADRN}. The velocities of the mobile components are defined by the interstitial velocity $v_\text{in}$ and their diffusion coefficients by the apparent diffusion coefficient $D_\text{app}$. The adsorption between the mobile and stationary phases can be compactly expressed by the stoichiometric form~\eqref{Rstochiometry}. The stoichiometric matrix is defined
\begin{align}
    \nu = \left[\begin{array}{cc}
    -\phi I; & I
    \end{array}\right],
\end{align}
where $\phi = (1-\varepsilon)/\varepsilon$ denotes the liquid volume fraction characterized by the column porosity $\varepsilon$. The associated kinetic rates $r = \partial_t q$ are given by a steric mass action model. For a protein $i \in \mathcal{C}_q$, it is defined by
\begin{align}
    k_{\text{kin},i} \partial_t q_i &=  k_{\text{eq},i}\bigg( \Lambda - \sum_{j \in \mathcal{C}_q}\left(\nu_j +\sigma_j\right) q_i \bigg)^{\nu_i}c_i \nonumber\\
    &\quad- q_i c_\text{NaCl}^{\nu_i},
\end{align}
where $\nu_i$ denotes the proteins characteristic charge, $\sigma_i$ its steric shielding factor, and $\Lambda$ the column's ionic capacity.
For the salt, it is defined by
\begin{align}
    \partial_tq_{\text{NaCl}} = - \sum_{j \in \mathcal{C}_q}\nu_j \partial_tq_j.
\end{align}